# The general solution of Schrödinger's differential equation

Nikos Bagis

**Aristotele University of Thessaloniki Greece**
**Department of Informatics**
**bagkis@hotmail.com**

**Abstract**
In this note we solve theoretically the Schrödinger's differential equation using results based on our previous work which concern semigroup operators. Our method does not use eigenvectors or eigenvalues and the solution depends only from the selected base of the Hilbert space.

**Introduction**

Recall that the Schrödinger equation describes the total energy of a particle in terms of potential and dynamical energy:

$$i\hbar \frac{\partial}{\partial t} \Psi = -\frac{\hbar^2}{2m} \nabla^2 \Psi + V(x)\Psi$$

where $\hbar$ is Planks constant, $m$ is the mass of particle, $V(x)$ is the potential, $\Psi(x,t) \in L^2(\mathbb{R}^n \times \mathbb{R}^+)$, $x$ is the position of the particle at time $t$. The evolution of the quantum system is expressed by $e^{-itH/\hbar}$, where

$$H = -\frac{\hbar^2}{2m} \nabla^2 + V(x)$$

is the Hamiltonian of the quantum system.

In this note we solve the equation for a particle moving on the real line, with an arbitrary potential $V(x)$, independent of time. The results "if they are correct" can be generalized and in higher dimensions.

**The general case of Schrödinger's equation**

**Definition 1.**
Let $f(.)$, $h_k(.)$, $G_k(.) \in L^2(\mathbb{R})$, $k \in \mathbb{Z}$. Then define $U_k$ the family of operators

$$(U_k f) = \langle f \mid h_k \rangle \int_{-\infty}^{\infty} G_k(x-t) f(t) dt$$

Let $T$ be an operator: $f \to Tf : L^2(\mathbb{R}) \to L^2(\mathbb{R})$ :

$$Tf(x) = \sum_{k=-\infty}^{\infty} c_k (U_k f)(x)$$

**Definition 2.**

We call "well posed" the class of all operators $S$ which are

$$f \to Sf : L^2(\mathbb{R}) \to L^2(\mathbb{R})$$

such that

$$Sf(x) = \int_{-\infty}^{\infty} f(t)K(t,x)dt$$

and where the kernel $K(t,x)$ has expansion in some base of $L^2(\mathbb{R})$, $\{e_k(x)\}_{k \in \mathbb{Z}}$, in the sense that

$$K(t,x) = \sum_{n=-\infty}^{\infty} Y_{k,l} e_k(t)e_l(x)$$

**Lemma.**

The operator $T$ is well posed.

**Proof.**

Using the isometric property of the Fourier Transform we have

$$Sf(x) = \int_{-\infty}^{\infty} f(t)\left(\sum_{k,l=-\infty}^{\infty} Y_{k,l} e_k(t)e_l(x)\right)dt$$

or

$$Sf^\wedge(\gamma) = \sum_{k,l=-\infty}^{\infty} Y_{k,l} \langle f \mid e_k \rangle e_l^\wedge(\gamma) = \sum_{k,l=-\infty}^{\infty} Y_{k,l} \langle f^\wedge \mid e_k^\wedge \rangle e_l^\wedge(\gamma)$$

But

$$Tf^\wedge(\gamma) = \sum_{k=-\infty}^{\infty} \frac{c_k}{2\pi} \langle f^\wedge \mid h_k^\wedge \rangle f^\wedge(\gamma)G_k^\wedge(\gamma)$$

We expand $h_k$ into orthonormal series of $e_k$ and finally we get

$$Tf^\wedge(\gamma) = f^\wedge(\gamma) \sum_{k,l=-\infty}^{\infty} \frac{c_l}{2\pi} \langle h_l^\wedge \mid e_k^\wedge \rangle \langle f^\wedge \mid e_k^\wedge \rangle G_l^\wedge(\gamma)$$

Setting $G_l(x) = e_l(x)$ and using the Riesz-Fischer Theorem ([R,N] pg. 70) for suitable $h_k$ we get

$$Tf^\wedge(\gamma) = f^\wedge(\gamma)Sf^\wedge(\gamma) : (1)$$

Hence

$$Tf(x) = \int_{-\infty}^{\infty} Sf(t)f(x-t)dt$$

Thus when we know the value of the one operator we can solve equivalently to find the other.

**Note.** More precisely in the proof of the Main Theorem we see how operators like $T$ are related to operators such $S$.

**Theorem. (Submitted in the Journal of Wavelets Theory and Applications)**

If $T$ is an operator as above then the solution of

$$\frac{\partial u(x,t)}{\partial t} = Tu(x,t) + g(x)$$

With initial condition $u(x,0) = f\left(x\right) \in L^2(\mathbb{R})$ is

$$u(x,t) = \frac{1}{2\pi} \int_{-\infty}^{\infty} f^{\wedge}(\gamma) e^{i\gamma x} e^{tK(\gamma)} d\gamma + \frac{1}{2\pi} \int_{0}^{t} \left( \int_{-\infty}^{\infty} g^{\wedge}(\gamma) e^{i\gamma x} e^{(t-s)K(\gamma)} d\gamma \right) ds \; : (2)$$

where $K(\gamma) = \sum_{k=-\infty}^{\infty} c_k h_k{}^{\wedge}(\gamma) G_k{}^{\wedge}(\gamma)$ and the operator $T$ is as in Definition 1.

**Main Theorem.**
The differential equation of Schrödinger's read as

$$\frac{\partial \Psi(x,t)}{\partial t} = -\frac{\partial^2 \Psi(x,t)}{\partial x^2} + V(x)\Psi(x,t) \; : (3)$$

If $\{e_k(x)\}_{k \in \mathbb{Z}}$ is an arbitrary base of $L^2(\mathbb{R})$ then the general solution of (3) is

$$\Psi(x,t) = \frac{1}{2\pi} \int_{-\infty}^{\infty} f^{\wedge}(\gamma) e^{i\gamma x} \exp\left( \sum_{k,m=-\infty}^{\infty} \left\langle \frac{e_k{}^{\wedge}(\gamma)}{e_m{}^{\wedge}(\gamma)} \middle| He_m{}^{\wedge}(\gamma) \right\rangle e_m{}^{\wedge}(\gamma) e_k{}^{\wedge}(\gamma) \right) d\gamma \; : (4)$$

Where the initial condition is $\Psi(x,0) = f^{\wedge}(\gamma)$.

**Proof.**
In view of Theorem 1 it is sufficient to write the Hamiltonian in the form

$$Sf(x) = \int_{-\infty}^{\infty} f(t) R(t,x) dt = Hf(x) = -\frac{d^2 f}{dx^2} + V(x)f(x) \; : (5)$$

and find the function $R(t,x)$. But as someone can see with

$$R(t,x) = -\delta''(x-t) + V(x)\delta(x-t)$$

where $\delta$ is the Dirac Delta function, we have

$$R_{k,l} = \int_{-\infty}^{\infty} \int_{-\infty}^{\infty} R(t,x) e_k(t) e_l(x) dt dx =$$

$$= \int_{-\infty}^{\infty} \int_{-\infty}^{\infty} -\delta''(x-t) e_k(t) dt e_l(x) dx + \int_{-\infty}^{\infty} \int_{-\infty}^{\infty} V(t)\delta(x-t) e_k(t) dt e_l(x) dx$$

$$= -\int_{-\infty}^{\infty} e_k''(x) e_l(x) dx + \int_{-\infty}^{\infty} V(x) e_k(x) e_l(x) dx \,.$$

Thus we see that $S$ is well posed and

$$R_{k,l} = \left\langle He_k \mid e_l \right\rangle \quad : (6)$$

where $H$ is the Hamiltonian of the Quantum System.

The only thing that lefts is how we can find the operator $T$, when we know $S$. For to construct the operator $T$ we know that

$$Tf^{\wedge}(\gamma) = f^{\wedge}(\gamma) \sum_{k,l=-\infty}^{\infty} \frac{c_l}{2\pi} \left\langle h_l{}^{\wedge} \mid e_k{}^{\wedge} \right\rangle \left\langle f^{\wedge} \mid e_k{}^{\wedge} \right\rangle G_l{}^{\wedge}(\gamma)$$

We choose $G_l(x) = e_l(x)$ and $T$ becomes:

$$Tf^\wedge(\gamma) = f^\wedge(\gamma) \sum_{k,l=-\infty}^{\infty} \frac{c_l}{2\pi} \langle h_l^\wedge \mid e_k^\wedge \rangle \langle f^\wedge \mid e_k^\wedge \rangle e_l^\wedge(\gamma)$$

We try to find the parameters of:

$$f^\wedge(\gamma) \left( \int_{-\infty}^{\infty} f(t) K(t,x) dt \right)^\wedge (\gamma) = Sf^\wedge(\gamma) \quad : \text{(a)}$$

and

$$f^\wedge(\gamma) \sum_{k,l=-\infty}^{\infty} \frac{c_l}{2\pi} \langle h_l^\wedge \mid e_k^\wedge \rangle \langle f^\wedge \mid e_k^\wedge \rangle e_l^\wedge(\gamma) = Sf^\wedge(\gamma) \quad : \text{(b)}$$

We start with (a), which can be rewritten as:

$$f^\wedge(\gamma) \sum_{k,l=-\infty}^{\infty} K_{k,l} \langle f^\wedge \mid e_k^\wedge \rangle e_l^\wedge(\gamma) = \sum_{k,l=-\infty}^{\infty} R_{k,l} \langle f^\wedge \mid e_k^\wedge \rangle e_l^\wedge(\gamma)$$

Set $f^\wedge(\gamma) = e_m^\wedge(\gamma)$, then

$$e_m^\wedge(\gamma) \sum_{k,l=-\infty}^{\infty} K_{k,l} \delta_{m,k} e_l^\wedge(\gamma) = \sum_{k,l=-\infty}^{\infty} R_{k,l} \delta_{m,k} e_l^\wedge(\gamma)$$

or

$$e_m^\wedge(\gamma) \sum_{l=-\infty}^{\infty} K_{m,l} e_l^\wedge(\gamma) = \sum_{l=-\infty}^{\infty} R_{m,l} e_l^\wedge(\gamma)$$

or if (i) : $e_m^\wedge(\gamma)$ has no real roots:

$$K_{m,s} = \sum_{l=-\infty}^{\infty} R_{m,l} \left\langle \frac{e_l^\wedge(\gamma)}{e_m^\wedge(\gamma)} \middle| e_s^\wedge(\gamma) \right\rangle = \sum_{l=-\infty}^{\infty} R_{m,l} \left\langle \frac{e_s^\wedge(\gamma)}{e_m^\wedge(\gamma)} \middle| e_l^\wedge(\gamma) \right\rangle \quad : (7)$$

(ii) : provided that the sum converges.
But

$$Tf^\wedge(\gamma) = f^\wedge(\gamma) \left( \int_{-\infty}^{\infty} f(t) K(t,x) dt \right)^\wedge (\gamma) = f^\wedge(\gamma) \sum_{m,s=-\infty}^{\infty} K_{m,s} \frac{1}{2\pi} \langle f^\wedge \mid e_m^\wedge \rangle e_s^\wedge(\gamma)$$

or from (a), (b), (7), (we seek for the $h_k$):

$$\sum_{m,s=-\infty}^{\infty} \frac{c_s}{2\pi} \langle h_s^\wedge \mid e_m^\wedge \rangle \langle f^\wedge \mid e_m^\wedge \rangle e_s^\wedge(\gamma) =$$

$$= \frac{1}{2\pi} \sum_{m,s=-\infty}^{\infty} \langle f^\wedge \mid e_m^\wedge \rangle \left( \sum_{l=-\infty}^{\infty} R_{m,l} \left\langle \frac{e_s^\wedge(\gamma)}{e_m^\wedge(\gamma)} \middle| e_l^\wedge(\gamma) \right\rangle \right) e_s^\wedge(\gamma)$$

Thus we choose

$$\langle h_s^\wedge \mid e_m^\wedge \rangle = \sum_{l=-\infty}^{\infty} R_{m,l} \left\langle \frac{e_s^\wedge(\gamma)}{e_m^\wedge(\gamma)} \middle| e_l^\wedge(\gamma) \right\rangle = \left\langle \frac{e_s^\wedge(\gamma)}{e_m^\wedge(\gamma)} \middle| \sum_{l=-\infty}^{\infty} R_{m,l} e_l^\wedge(\gamma) \right\rangle =$$

$$= \left\langle \frac{e_s^\wedge(\gamma)}{e_m^\wedge(\gamma)} \middle| \sum_{l=-\infty}^{\infty} \langle He_m \mid e_l \rangle e_l^\wedge(\gamma) \right\rangle = \left\langle \frac{e_s^\wedge(\gamma)}{e_m^\wedge(\gamma)} \middle| He_m^\wedge(\gamma) \right\rangle$$

Note in the general case we have $\langle h_s^\wedge \mid e_m^\wedge \rangle = \left\langle \frac{e_s^\wedge(\gamma)}{e_m^\wedge(\gamma)} \middle| He_m^\wedge(\gamma) \right\rangle = \mu_{s,m}$

Hence clearly

$$h_s\,{}^\wedge(\gamma) = \sum_{m=-\infty}^{\infty} \mu_{m,s} e_m\,{}^\wedge(\gamma) = \sum_{m=-\infty}^{\infty} \left\langle \frac{e_s\,{}^\wedge(\gamma)}{e_m\,{}^\wedge(\gamma)} \middle| He_m\,{}^\wedge(\gamma) \right\rangle e_m\,{}^\wedge(\gamma)$$